\documentclass[11pt]{article}
\usepackage{amsmath,amsthm,amsfonts}
\usepackage{color}
\usepackage{ amssymb }

\numberwithin{equation}{section} \allowdisplaybreaks

\newcommand{\N}{\mathbb N}

\newtheorem{theorem}{\sc Theorem}[section]
\newtheorem{lemma}[theorem]{\sc Lemma}
\newtheorem{proposition}[theorem]{\sc Proposition}
\newtheorem{corollary}[theorem]{\sc Corollary}
\newtheorem{definition}[theorem]{\sc Definition}
\newtheorem{example}[theorem]{\sc Example}
\newtheorem{remark}[theorem]{\sc Remark}
\newtheorem{remarks}[theorem]{\sc Remarks}
\newtheorem{examples}[theorem]{\sc  Examples}
\newtheorem{Application}[theorem]{\sc  Application}
\newcommand{\bremarks}{\begin{remarks}\rm}
\newcommand{\eremarks}{\end{remarks}}
\newcommand{\bexamples}{\begin{examples}\rm}
\newcommand{\eexamples}{\end{examples}}
\newcommand{\bapp}{\begin{Application}\rm}
\newcommand{\eapp}{\end{Application}}

\def \n {\noindent}
\def \n {\noindent}

\newcommand{\bet}{\begin{theorem}}
\newcommand{\eet}{\end{theorem}}
\newcommand{\blm}{\begin{lemma}}
\newcommand{\elm}{\end{lemma}}
\newcommand{\bprop}{\begin{proposition}}
\newcommand{\eprop}{\end{proposition}}
\newcommand{\bcor}{\begin{corollary}}
\newcommand{\ecor}{\end{corollary}}
\newcommand{\bdf}{\begin{definition}\rm}
\newcommand{\edf}{\end{definition}}
\newcommand{\bp}{\begin{proof}}
\newcommand{\ep}{\end{proof}}
\newcommand{\bex}{\begin{example}\rm}
\newcommand{\eex}{\end{example}}
\newcommand{\bremark}{\begin{remark}\rm}
\newcommand{\eremark}{\end{remark}}

\newcommand{\nul}[1]{\mathrm{N}( {#1} )}
\newcommand{\ran}[1]{\mathrm{R}({#1})}

\newcommand{\codim}{\mathrm{codim \,}}

\newcommand{\ch}{{ \mathcal C}(X)}

\begin{document}

\title {On the B-discrete spectrum }

\author { M. Berkani}

\date{}

\maketitle

\begin{abstract}

\noindent

 In this paper, we introduce  the B-discrete spectrum of an
unbounded closed operator and we prove that a closed operator has
a purely B-discrete spectrum if and only if it has a meromorphic
resolvent.  After that, we study the stability of the B-discrete
spectrum under several type of perturbations and we establish that
 two closed invertible linear operators having quasisimilar
totally paranormal inverses have equal spectra and B-discrete
spectra.

\end{abstract}

\renewcommand{\thefootnote}{}

\footnotetext{\hspace{-7pt}2010 {\em Mathematics Subject
Classification\/}:  primary 47A10, 47A53.
\baselineskip=18pt\newline\indent {\em Key words and phrases\/}:
   B-discrete,  meromorphic, perturbation,  hereditarily normaloid }

\section{Introduction}

   Let $\ch$ be  the  set of  linear closed operators defined
from a Banach space  $X$ to $X$ and $L(X)$ be the Banach algebra
of  bounded linear  operators defined from   $X$ to $X.$ We write
${\mathrm D}(T)$, $\nul{T}$ and $\ran{T}$ for the domain,
nullspace and range of an operator $T\in \ch$. An operator $T\in
\ch$  is called a {\it Fredholm}  operator \cite{TLY} if both the
nullity  $n(T)=\dim \nul{T}$ of $T$ and the defect $d(T)=\codim
\ran{T}$  of $T$  are finite. The index $i(T)$ of a Fredholm
operator $T$ is defined by $i(T)=n(T)-d(T)$. It is well known that
if $T$ is a Fredholm operator, then $\ran{T}$ is closed.

The class of bounded linear  B-Fredholm operators, which is a
natural extension of the class of Fredholm operators  was
introduced in \cite{P7}, and the class of unbounded linear closed
B-Fredholm operators acting on a Banach space was studied in
\cite{P33}.

Recall \cite{CA-1} that a bounded linear  operator $T$  is called
a meromorphic operator if $ \lambda = 0$ is the only possible
point of accumulation of its spectrum $\sigma(T)$ and every
non-zero isolated point of $\sigma(T)$ is a pole of the resolvent
$R_\texttt{{\tiny $\mu $}}(T)= (T-\mu I)^{-1}$  of $T,$ which
is defined on the resolvent set $\rho(T)$ of $T.$  If we also
require that each non-zero eigenvalue of T has finite
multiplicity, then $T$ will be called a Riesz operator.

A first result linking bounded B-Fredholm operators to the class $
\mathfrak{M}$ of linear bounded meromorphics operators comes from
the following theorem, established in \cite[Theorem 2.11]{P13}.

\bet \label{thm1.1} Let $ T \in L(X).$  Then    $T$ is a
meromorphic operator if and
 only if  $ \sigma_{BF} (T)= \{ \lambda \in \mathbb{C} \mid  T - \lambda I \text {\,
 is not a B-Fredholm operator } \} \subset \{0\}.$
\eet

Recall that for $T\in \ch$,  its  {\it descent} $\delta(T)$ and
its  {\it ascent} $ a(T) $  are defined by
 $ \delta(T)=\inf  \{n\in \mathbb{N}: R(T^n) = R(T^{n+1})\}$
and
 $ a(T)= \inf \{n\in \mathbb{N} : N(T^{n})=N(T^{n+1})\}$. We set formally $\inf\emptyset =\infty$.

A closed linear operator $T\in \ch$ is said to be {\it Drazin
invertible} if $a(T)$ and $ \delta(T)$ are both finite.  In this
case and  if the resolvent set $\rho(T)$ of $T$  is nonempty, then
$a(T)= \delta(T),  R(T^{a(T)})$ is closed and $ X= R(T^{a(T)})
\oplus  N(T^{a(T)}).$

The {\it Drazin spectrum} of $T$ is defined by: $ \sigma_{\cal
D}(T)= \{\lambda\in\mathbb{C}: T-\lambda
   I \text{ not Drazin invertible} \}.$

The set  of  Browder operators is defined by $\mathcal{B}(X)=\{ T
\in \Phi (X) \mid  a (T)<\infty \,\, \text {and} \,\,\delta(T) <
\infty \}$ and
  the {\it  Browder spectrum} of $T$ is defined by: $ \sigma_{{\cal B}}(T)= \{\lambda\in\mathbb{C}:T-\lambda
   I\notin \mathcal{B}(X)\}.$

 For a closed linear operator $T \in \ch,$ the discrete spectrum
$\sigma_{ d}(T)$ of $T$ is defined as the set of all complex
numbers $\lambda$ in $\sigma(T)$ such that $ T- \lambda I $ is  a
Browder operator, that is $\sigma_{d}(T) = \sigma(T) \setminus
 \sigma_{{\cal B}}(T),$ the complement of the Browder spectrum in
 the spectrum.

Analogously,  we define here the  B-discrete spectrum  for closed
operators as a natural extension of the  discrete spectrum.

\bdf \label{def1} Let $T \in \ch $.
 Then the   B-discrete
spectrum $\sigma_{bd}(T)$ of $T$  is defined by $\sigma_{bd}(T)=
\sigma(T) \setminus \sigma_{\cal D}(T),$  the complement of the
Drazin spectrum in
 the spectrum.\edf

It's clear that for $T \in \ch,$ its discrete spectrum
$\sigma_{{d}}(T)$  is a subset of its  B-Discrete spectrum
$\sigma_{bd}(T).$

\bdf \label{def9} We will say that $T$ has a purely B-discrete
spectrum if $\sigma(T)= \sigma_{bd}(T), $ and  that $T$ has a
purely discrete spectrum if $\sigma(T)= \sigma_{d}(T).$ \edf

\bex \label{ex1}
 An illustrating example of an operator with purely B-discrete
spectrum, is given by the Schr$\ddot{o}$dinger operator with a
constant magnetic field $B\neq 0 $ in $\mathbb{R}^2$ defined by
$S_{B}= (\frac{1}{i}\frac{\partial}{\partial
x_1}-\frac{Bx_2}{2})^2 +(\frac{1}{i}\frac{\partial}{\partial
x_2}+\frac{Bx_1}{2})^2. $ Then  from \cite[Example 4, p.134]{BH}
$\sigma(S_B)= \{ (2k+1)\mid B \mid \ | k \in \mathbb{N} \}.$ Its
B-discrete spectrum coincides exactly with the set of  its Landau
levels,  (\cite[ p.136]{BH}) while its discrete spectrum is the
empty set. Moreover, each eigenvalue of $S_B$ has an infinite
multiplicity. As $S_B$ is self-adjoint, we have $\sigma(S_B)=
\sigma_{bd}(S_B),$ (See Corollary \ref{NN}). Thus $L$ has a purely
B-discrete spectrum, but its discrete spectrum is empty.
 \eex

  The discrete spectrum has important applications in the study
of physical operators. However, as shown by Example \ref{ex1}, the
discrete spectrum may be empty while we have discrete energy
levels! Moreover, as shown by Example \ref{Qnil}, the discrete
spectrum does not give a clear idea on the nature of isolated
points of the spectrum: are there poles of infinite rank or
essential singularities of the resolvent?  But an isolated point
of the spectrum which is not in the B-discrete spectrum is always
an essential singularity.

\bdf \label{def6}  Let $T \in \ch, $  with a non-empty resolvent
set. We will say that $T$ has
 a meromorphic (resp. Riesz or  compact) resolvent if there exists a
scalar $\lambda $ in the resolvent set $ \rho(T)$  of  $T$  such
that $ (T-\lambda I)^{-1}$ is a bounded linear  meromorphic (resp.
Riesz or compact) operator. \edf

 \bremark It's easily seen that if $T$ has a meromorphic (resp. Riesz or compact)  resolvent, then
 for all
scalar $\lambda $ in the resolvent set $ \rho(T)$  of  $T,$  $
(T-\lambda I)^{-1}$ is a  bounded linear meromorphic (resp. Riesz
or compact) operator.

 \eremark

  In the second section of this paper, we characterize closed invertible operators with non-empty resolvent set
having a purely B-discrete spectrum, by showing that this the case
if and only if the operator considered has a  meromorphic
resolvent and  if and only if its B-Fredholm spectrum $
\sigma_{BF} (T)$ is empty. We show also that if  $ T \in \ch$ has
a nonempty resolvent set and  $ \lambda$ is an isolated point of
its spectrum, then $ \lambda $ is in its B-discrete spectrum if
and only if $ T-\lambda I$  is a B-Fredholm operator. When $T$ is
an hereditarily normaloid operator  (Definition \ref
{h-normaloid}), then the B-discrete spectrum $\sigma_{ bd}(T)$ of
$T$ is the set of all isolated points of its spectrum $\sigma(T).$

  In the third section, we study the stability of  the
 B-discrete spectrum under several type of perturbations. As an example   of the
results obtained, we show that if  $A$ and $T$ are two commuting
closed linear operators with nonempty resolvent sets, and if for
some $\lambda\in\rho(A)\cap\rho(T)$ the operator $(\lambda
I-A)^{-1}-(\lambda I-T)^{-1} $  is of finite rank,  then $A$ has a
purely B-discrete spectrum if and only if $T$ has a purely
B-discrete spectrum. Moreover, we prove two spectral mapping
theorems for the B-discrete spectrum.

 In the fourth section, we show that if  $S, T$ are two closed
invertible linear operators having quasisimilar totally paranormal
 inverses (Definition \ref{paranormal}), then their spectra and
their B-discrete spectra are equal.

\section{B-discrete spectrum}

We begin this section by characterizing operators with purely
B-discrete spectrum.

\bet \label{thm10} Let $ T \in \ch$ with a nonempty resolvent set.
Then $T$ has a purely B-discrete spectrum  if and only if $T$ has
a  meromorphic resolvent. \eet

\bp  Suppose that $T$ has a purely B-discrete spectrum.  So for
all $ \lambda \in \mathbb{C}, \,  T - \lambda I $ is Drazin
invertible operator. Since the resolvent set of $T$ is non-empty,
there exists $\mu \in \mathbb{C},$ such that $T- \mu I$ is
invertible. From \cite[Theorem  3.6]{P40}, $(T- \mu I)^{-1}
-\frac{1}{\lambda} I$ is  Drazin invertible for all $\lambda \neq
0,$ and $\frac{1}{\lambda}$  is a pole of the resolvent of $(T-
\mu I)^{-1}.$ So $(T- \mu I)^{-1}$ is a meromorphic operator.
Hence $T$ has  a meromorphic resolvent.

 Conversely,  if
$T$ has  a meromorphic resolvent, we can assume without loss of
generality that
  $ T$ is invertible and that $ T^{-1}$ is a meromorphic operator. If
  $\lambda \notin \sigma(T),$ then $T -\lambda I$ is invertible.  If
  $\lambda \in \sigma(T),$ then  $ \lambda \neq 0.$ Since $T^{-1}$ is a meromorphic operator,
  then from \cite[Theorem  3.6]{P40}, $\frac{1}{\lambda }$ is a pole  of $T^{-1}$
  and again from \cite[Theorem  3.6]{P40},
    $\lambda$ is a pole of $T.$  Therefore $ T - \lambda I $ is Drazin invertible   for all $
 \lambda \in \mathbb{C}$ and $T$ has a purely B-discrete spectrum.
\ep

\bremark Analogous result of Theorem \ref{thm10} for the discrete
spectrum can be deduced from  \cite[Theorem 2]{KLY}, that  is $T$
has a purely discrete spectrum  if and only if $T$ has a Riesz
resolvent.

\eremark

\bcor Let $ T \in \ch $ with a nonempty resolvent set. Then $T$
has a purely B-discrete spectrum if and only if $ \sigma_{BF} (T)=
\emptyset.$  \ecor

\bp   Suppose that $T$ has a purely B-discrete spectrum.  Then for
all $ \lambda \in \mathbb{C}, \,  T - \lambda I $ is Drazin
invertible. From \cite[Theorem 2.9]{P33}, it follows that  $T
-\lambda I$ is a B-Fredholm operator. Hence  $ \sigma_{BF} (T)=
\emptyset.$

Conversely assume that $ \sigma_{BF} (T)= \emptyset,$ and let
$\mu$ be in the resolvent set of $T.$  Then from \cite[Theorem
3.6]{P40}, $ \sigma_{BF}( (T- \mu I)^{-1} ) \subset \{0\}$ and
from Theorem \ref{thm1.1}, $(T- \mu I)^{-1}$ is a meromorphic
operator. Hence $T$ has  a meromorphic resolvent. From Theorem
\ref{thm10}, it follows that $T$ has purely B-discrete spectrum.

\ep

\begin{example} Let $A$  be the shift operators on the Hilbert space $l^{2}(\mathbb{N})$.
\begin{eqnarray*}
  A:D(A)\subset l^{2}(\mathbb{N}) &\longrightarrow & l^{2}(\mathbb{N}) \\
    x=(x_{n})_{n\geq 0} &\longmapsto & Ax=(0,x_{1},2x_{2},3x_{3},.......),
\end{eqnarray*}
where  $D(A)=\{x=(x_{n})_{n\geq 0}\in l^{2}(\mathbb{N}):
\displaystyle\sum_{n\geq 0}n^{2}|x_{n}|^{2}<\infty\}$

It is easy to verify that the complex number $i$ $(i^{2}=-1)$ is
in $\rho(A)$. Moreover one  can check easily that $\sigma(A)=\N $
and that $\sigma_{BF}(A)=\emptyset.$
 Hence $A$ has a purely B-discrete spectrum.
\end{example}

\bcor Let $ T \in \ch,$ with a nonempty resolvent set. If $T$ has
a purely discrete spectrum, then $T$  has a purely B-discrete
spectrum and $ \sigma_{bd}(T)= \sigma_{d}(T)$. \ecor

\bp It follows from \cite[Theorem 2]{KLY}, that if $T$ has purely
discrete spectrum, then $T$ has a Riesz resolvent. So $T$ has a
meromorphic resolvent. Thus  $T$  has a purely B-discrete
spectrum. In this case, we have $  \sigma(T)= \sigma_{ d}(T),$ so
$\sigma_{ B}(T)= \emptyset.$ Hence $\sigma_{D}(T)= \emptyset$
 and $
\sigma(T)= \sigma_{bd}(T). $ \ep

\noindent As shown by the example \ref{ex1}, the converse of the
previous corollary is not true in general.

\vspace{2mm}

The following example shows that even the spectrum  of an operator
is discrete, its B-discrete spectrum ( and also its discrete
spectrum) could be empty.

\bex\label{Qnil} Let $Q$ be defined for each $x=(\xi_i)\in\ell^1$
by
\begin{equation*}
Q(\xi_1,\xi_2,\xi_3,\dots, \xi_k,\dots) = ( 0, \alpha_1\xi_1,
\alpha_2\xi_2,\dots, \alpha_{k-1}\xi_{k-1},\dots),
\end{equation*}
where $(\alpha_i)$ is a sequence of complex numbers such that
$0<|\alpha_i|\le 1$ and $\sum_{i=1}^\infty |\alpha_i|<\infty$.  We
observe that
\begin{equation*}
\overline{R(Q^n)} \ne R(Q^n), \quad n=1,2,\dots
\end{equation*}
Indeed, for a given $n\in\N$ let $x^{(n)}_k = (1,\dots ,
1,0,0,\dots)$ (with $n+k$ times $1$).  Then the limit
$y^{(n)}=\lim_{k\to\infty}Q^nx^{(n)}_k$ exists and lies in
$\overline{R(Q^n)}$.  However, there is no element
$x^{(n)}\in\ell^1$ satisfying the equation $Q^nx^{(n)}=y^{(n)}$ as
the algebraic solution to this equation is  $(1,1,1,\dots
)\notin\ell^1$.

It is easy to see that  $ \sigma (Q)= \sigma_{D} (Q)= \{ 0\}$ and
so $\sigma_{d}(T)=\sigma_{bd}(T)= \emptyset. $

 \eex

In the next theorem, we give a necessary and sufficient condition
under which an isolated point in the spectrum of a closed operator
is in its B-discrete spectrum.

\bet \label{iso} Let $ T \in \ch,$ with a nonempty resolvent set
and let $ \lambda$ be an isolated  point of its spectrum. Then $
\lambda $ is in its B-discrete spectrum if and only if $ T-\lambda
I$  is a B-Fredholm operator. \eet

\bp If $ \lambda \in \sigma_{bd}(T),$ then $ \lambda $ is a pole
of its resolvent.  From \cite[Theorem 1.2]{LY}, there exists an
integer $p$ such that  $ R((T - \lambda I)^p)$  is closed and $X=
N((T-\lambda I)^p)\oplus R((T - \lambda I)^p).$ From \cite[Theorem
2.4]{P33}, it follows that  $ T-\lambda I$ is a B-Fredholm
operator.

Conversely assume that $ \lambda $ is isolated  in $\sigma(T)$ and
that $ T-\lambda I$ is a B-Fredholm operator. We can assume
without loss of generality that $T$ is invertible, thus $ \lambda
\neq 0 $ and from \cite[Theorem 3.6]{P40}, $T^{-1} -
\frac{1}{\lambda} I $ is a B-Fredholm operator. From \cite[Theorem
2.3]{P13}, it follows that $\frac{1}{\lambda}$ is a pole of
$T^{-1}$ and from \cite[Theorem 3.6]{P40} it follows that $\lambda
$ is a pole of $T.$  So  $ \lambda $ is in the B-discrete spectrum
of $T.$ \ep

\bdf \label{h-normaloid}

 A bounded linear operator $T$ acting on a Banach space is said to be
normaloid if $\rho(T) =\| T\|$, where $\rho(T)$ is the spectral
radius of $T$ or equivalently $\| T^n\| = \| T\| ^n $ for every $
n \in \mathbb{N}.$ We say that $T$ is hereditarily  normaloid if
the restriction of $T$ to any invariant subspace under $T$  is
again normaloid.

\edf


\bdf \label{paranormal}

A bounded operator $T$ on a Banach space is said to be paranormal
if $\| Tx\|^2 \leq \| T^2 x \|\| x\|$ for all $x \in X . $ $T$ is
said to be totally paranormal if $T- \lambda I$ is paranormal for
all $ \lambda \in \mathbb{C}.$   \edf

It is easy to see that every paranormal operator $T$ is normaloid.
Since the restriction of a paranormal operator $T$ to any closed
invariant subspace under $T$ is also paranormal, then a paranormal
operator is hereditarily normaloid.

A good example of totally paranormal operators is given by the
class of hyponormal operators. Recall that a bounded linear
operator $T$ acting on Hilbert space $H$ is said to be hyponormal
if $\| T^*x \| \leq  \| Tx \|,$  for all $x \in H, $ where $ T^*$
is the adjoint of $T.$ From \cite{FUR}, we know that Hyponormal
$\subsetneq $  Paranormal $\subsetneq$ hereditarily normaloid
$\subsetneq$ normaloid and the inclusions are all proper.

In the case of a closed  invertible operator with a   normaloid
inverse, we can determine precisely its B-discrete spectrum.

\bet \label{HN}Let $ T \in \ch$ be a closed  invertible operator
having an hereditarily normaloid inverse. Then the B-discrete
spectrum $\sigma_{ bd}(T)$  of $T$ is the set of all isolated
points of the spectrum $\sigma(T)$ of $T.$ \eet

\bp  If $ \lambda \in \sigma_{bd}(T),$ then $ \lambda $ is a pole
of its resolvent and so $\lambda $ is an isolated point of its
spectrum. Conversely assume that  $ \lambda$ is  an isolated point
of the spectrum of $T.$  As $T^{-1}$ is hereditarily normaloid,
then from \cite [Lemma 2.1]{DD} $\frac{1}{\lambda}$  is a pole of
$T^{-1}$ and from From \cite[Theorem 3.6]{P40}$ \lambda$ is a pole
of $T.$ So $ \lambda$ is in $\sigma_{bd}(T).$

\ep

\bcor \label{NN} Let $T$ be a normal operator  acting on a Hilbert
space $H.$
 Then the B-discrete spectrum
$\sigma_{ bd}(T)$  of $T$ is the set of all isolated points of the
spectrum $\sigma(T)$ of $T.$ \ecor

\bp
   Since the resolvent set $\rho(T)$ of
$T$ is nonempty, we can assume without loss of generality that $
T$ is invertible. Since $T$ is normal, then from \cite[Theorem
5.42]{Wei}, $ T^{-1} $ is a bounded   normal operator. Since a
normal operator is normaloid, then from Theorem \ref{HN},
$\sigma_{ bd}(T)$ is the set of all isolated points of the
spectrum $\sigma(T)$ of $T.$

\ep

\vspace{5mm}

\noindent Corollary  \ref{NN}  applies  in particular to the case
of self-adjoint operators and in particular to self-adjoint
Schr{\"{o}}dinger operators \cite{Wei}.

 \section{Stability of the B-discrete spectrum }

\hskip0.7cm In this Section, we study  the stability  of the
B-discrete spectrum   of a closed linear operator acting on a
Banach space $X,$ under the effect of several type of
perturbations

\begin{definition}\cite{MJT}\label{0.0.0}
Let $X$ be a Banach space, $A: D(A)\subset X\longrightarrow X$
and\\ $T: D(T)\subset X\longrightarrow X$ two linear operators. We
say that $A$ commutes with $T$ and we denote $AT=TA$, if\\$(i)$
$D(A)\subset D(T)$.\\$(ii)$ $Tx\in D(A)$ whenever $x\in
D(A)$.\\$(iii)$ $AT=TA$ on $\{x\in D(A),~~Ax\in D(T) \}$.

\end{definition}

It is easily seen that if  $A$ and $T$ are two  commuting closed
linear operators on a Banach space $X$ and if
$\lambda\in\rho(A)\cap\rho(T),$ then $(\lambda I-A)^{-1}(\lambda
I-T)^{-1}= (\lambda I-T)^{-1}(\lambda I-A)^{-1}$.

\bet  \label{nil} Let $A$ and $T$ be two closed linear operators
 with a nonempty resolvent sets. If $AT=TA$ and
$(\lambda I-A)^{-1}-(\lambda I-T)^{-1}$ is a nilpotent operator
for some $\lambda\in\rho(A)\cap\rho(T)$, then
$\sigma_{bd}(A)=\sigma_{bd}(T).$

\eet

\bp Without loss of generality, we can assume   that $\lambda=0$.
Let $\mu\in\mathbb{C}\backslash\{0\}$. Since $T^{-1}- A^{-1}$ is a
nilpotent operator commuting with $ A^{-1},$ then $\sigma(A^{-1} )
= \sigma(A^{-1}+ (T^{-1}- A^{-1}))= \sigma(T^{-1}).$ Moreover and
from \cite[Theorem 3.2] {P22}  we know that $\sigma_{D}(A^{-1})=
\sigma_{D}(A^{-1}+ (T^{-1}- A^{-1}))= \sigma_{D}(T^{-1}).$ Thus
$\sigma(A)= \sigma(T)$ and  $\sigma_D(A)= \sigma_D(T).$
Consequently we have $\sigma_{bd}(A)=\sigma_{bd}(T).$ \ep

\begin{corollary}\label{nil1}
Let $A\in C(X)$ with a nonempty resolvent set  and let $Q\in L(X)$
be a nilpotent operator satisfying $AQ=QA$. Then
$\sigma_{bd}(A+Q)=\sigma_{bd}(A).$
\end{corollary}
\n {\bf Proof.}  Since $AQ=QA$ and $Q$ is a nilpotent operator,
then $(\mu I-A)^{-1}Q$ is also a nilpotent operator for all
$\mu\in\rho(A)$. Then  $I-(\mu I-A)^{-1}Q$ is an invertible
operator and so $\mu I-A-Q$ is also invertible. Moreover $(\mu
I-A-Q)^{-1}=(\mu I-A)^{-1}(I-(\mu I-A)^{-1}Q)^{-1}=(\mu
I-A)^{-1}\displaystyle\sum_{k=0}^{n}((\mu I-A)^{-1}Q)^{k}= (\mu
I-A)^{-1}+(\mu I-A)^{-1}Q\displaystyle\sum_{k=1}^{n-1}((\mu
I-A)^{-1})^{k}Q^{k-1}$ where $n$ is the nilpotent-index  of $(\mu
I-A)^{-1}Q$. Hence, $(\mu I-A-Q)^{-1}-(\mu I-A)^{-1}$ is
nilpotent. From Theorem \ref{nil}, we deduce that
$\sigma_{bd}(A+Q)=\sigma_{bd}(A)$.

\bet \label{1} Let $A$ and $T$ be two closed  linear operators on
a Banach space $X.$  If $AT=TA$ and  for some
$\lambda\in\rho(A)\cap\rho(T)$ the operator $(\lambda
I-A)^{-1}-(\lambda I-T)^{-1}$ has a power of finite rank, then $A$
has a purely  B-discrete spectrum if and only if $T$ has a purely
B-discrete spectrum.

\eet

\bp Without loss of generality we can assume that $\lambda=0$. Let
$\mu\in\mathbb{C}\backslash\{0\}$.  From \cite[Theorem 3.6]{P40},
we see that $\mu I-A$ is a B-Fredholm operator if and only if
$\mu^{-1}I-A^{-1}$ is
 a B-Fredholm   operator. Since $A^{-1}-T^{-1}$ has a power of
finite rank, then from \cite[Theorem 2.11]{ZJZ},
$\mu^{-1}I-A^{-1}$ is a
 B-Fredholm  operator
 if and only if $\mu^{-1}I-T^{-1}$ is a
 B-Fredholm operator.  From \cite[Theorem 3.6]{P40}  this is true
if and only if $\mu I-T$ is a B-Fredholm  operator. This shows
that
  $\sigma_{BF}(A)=\sigma_{BF}(T).$  Therefore
$A$ has a purely  B-discrete spectrum if and only if $T$ has a
purely  B-discrete spectrum. \ep

\vspace{2mm} We conclude this section with two spectral mapping
theorems for the B-discrete spectrum.

\bet Let $A\in C(X)$ be a densely defined closed operator such
that $\rho(A)\neq\emptyset$. Let $P(\lambda)$ be a polynomial with
complex coefficients. Then $\sigma_{BF}(P(A))= P(\sigma_{BF}(A)),$
and $A$ has a purely B-discrete spectrum if and only if $P(A)$ has
a purely B-discrete spectrum. \eet

\bp P(A) is well defined and is a closed operator. From
\cite[Theorem 3.2]{P44} we have $\sigma_{BF}(P(A))=
P(\sigma_{BF}(A)).$ Thus $\sigma_{BF}(A)= \emptyset$ if and only
if $\sigma_{BF}(P(A))= \emptyset.$

\ep

\vspace{1mm}

For an unbounded closed operator $A$ with non-empty resolvent set,
and a complex-valued functions $f$ holomorphic on an open set
containing $ \sigma(A) \cup \{\infty\},$ $f(A)$ may be defined by
the operational calculus introduced by Taylor in \cite{TLY}.

\bet  Let $A\in C(X)$ be a closed operator with non-empty
resolvent set and let $f$ be  complex-valued function holomorphic
on an open set containing $ \sigma(A) \cup \{\infty\}.$  If $f$ is
an univalent function, then  $\sigma_{bd}(f(A))= f(
\sigma_{bd}(A)).$ \eet

\bp Let $\widetilde\sigma(A)=\sigma(A)\cup\{\infty\}$ and
$\widetilde{\sigma_{D}}(A)=\sigma_{D}(A)\cup\{\infty\}$  be the
extended spectrum and  the extended Drazin spectrum of $A.$ From
\cite[Theorem 7]{GL}, we know that  $ f( \widetilde\sigma(A))=
\sigma(f(A)) $ and from  \cite[Theorem 4.1]{P43} $ f(
\widetilde{\sigma_{D}}(A))= \sigma_{D}(f(A)). $ Therefore we have
$\sigma_{bd}(f(A))= \sigma(f(A)) \setminus \sigma_{D}(f(A))= f(
\widetilde\sigma(A)) \setminus  f( \widetilde{\sigma_{D}}(A))= f(
\widetilde\sigma(A) \setminus \widetilde{\sigma_{D}}(A))=  f(
\sigma(A)) \setminus \sigma_{D}(A)),$ because $f$ is univalent.
Hence $\sigma_{bd}(f(A))= f( \sigma_{bd}(A)).$

 \ep

 \section{ B-discrete spectrum and quasi-similarity }

A bounded linear operator $ A : X \rightarrow Y$  from the Banach
spaces $ X$ to the Banach space $Y$ is said to be quasi-invertible
if it is injective and has dense range. Two bounded linear
operators $T \in L(X)$ and $S \in L(Y )$ are quasisimilar if there
exists quasi-invertible operators $A : X \rightarrow Y$ and $B : Y
\rightarrow  X$ such that $AT = SA$ and $BS = T B.$

As mentioned in  \cite[p.89]{CLA}, the same proof of \cite[Theorem
1]{CLA} proved for hyponormal operators holds also for   totally
paranormal operators. Thus we formulate the following result
without proof and we refer the reader to \cite{CLA}.

\bet \label{para}  If two bounded linear operators $T \in L(X)$
and $S \in L(Y )$ are totally paranormal and quasisimilar, then
they have the same spectrum. This is in particular the case of two
quasisimilar hyponormal operators. \eet

Using this result, we prove now the equality of the B-discrete
spectrum of two  quasi-similar operators.

\bet \label{bd-para} Let $S$ and $T$ be two  totally paranormal
and quasisimilar bounded linear operators acting on a Banach space
$X,$  then $\sigma_{bd}(S)= \sigma_{bd}(T).$ This is in particular
the case of two quasisimilar hyponormal operators.

\eet

\bp Since   $S$ and $T$  are quasisimilar, there exists
quasi-invertible operators $A : X \rightarrow Y$  and $B : Y
\rightarrow X $ such that $AT = SA$ and $BS = T B.$  In this case
the operators $ \tilde{S}_ {[n]} : \overline{R(S^ n)} \rightarrow
\overline{R(S^ n)}$ and $ \tilde{T}_ {[n] }: \overline{R(T^ n)}
\rightarrow \overline{R(T^ n)}$  defined as the restrictions of $
S$ and $T$ respectively to the closure of the ranges $R(S^n)$ and
$R(T^n)$, are also quasisimilar. Indeed, if we consider $ P =
A/_{\overline{R(T^n)}}: \overline{R(T^n)} \rightarrow \overline{
R(S^n)}$ and $ Q= B/_{\overline{R(S^n)}}: \overline{R(S^n)}
\rightarrow \overline{ R(T^n)}$, then it is easily seen that
$\overline{P(\overline{R(T^n)})} = \overline{R(S^n)},
\overline{Q(\overline{R(S^n)})} = \overline{R(T^n)}, P$ and $Q$
are both injective, $P \tilde{T}_ {[n] }= \tilde{S}_ {[n] } P$ and
$Q\tilde{S}_ {[n] }= \tilde{T}_ {[n] } Q.$

Let $\alpha \in \sigma_{bd}(T)$  be arbitrary. Then $T - \alpha I$
is Drazin invertible and $ a(T - \alpha I) = \delta(T - \alpha I)
= n < \infty. $  Since we are dealing with totally paranormal
operators,   we may assume without loss of generality that $\alpha
= 0.$ Therefore $R(T^n)$ is closed, and  $ \tilde{T}_ {[n] } :
R(T^n) \rightarrow R(T^n)$ is invertible. On the other hand,
$\tilde{T}_ {[n] }  : R(T^ n) \rightarrow R(T^n)$  and $
\tilde{S}_ {[n]} : \overline{R(S^ n)} \rightarrow \overline{R(S^
n)} $ are totally paranormal  quasisimilar operators. Since
$\tilde{T}_ {[n] }$ is invertible, then from Theorem \ref{para},
$\tilde{S}_ {[n]}$ is invertible. So $(\tilde{S}_ {[n]})^ n$ is
invertible and $ R((\tilde{S}_ {[n]})^ n) = \overline{R(S^n)}.$ As
$ R((\tilde{S}_ {[n]})^ n) \subset R(S^n),$ then $
\overline{R(S^n)} = R(S^n)$ and  $R(S^ n)$ is closed. Therefore $0
\notin \sigma_D(S).$  As we know from Theorem \ref{para} that
$\sigma(T) = \sigma_(S),$  then $0 \in \sigma_{bd}(S).$ Similarly,
we have $\sigma_{bd}(S) \subset \sigma_{bd}(T)
$ and so
$\sigma_{bd}(S) = \sigma_{bd}(T).$

\ep

\bet  Let $S, T$ be  two closed invertible linear operators having
quasisimilar totally paranormal inverses, then $\sigma(S)=
\sigma(T)$ and  $\sigma_{bd}(S)= \sigma_{bd}(T).$ \eet

\bp  Let $ U= S^{-1},$ and $V= T^{-1},$ then $U$ and $V$ are
totally paranormal and quasi-similar. From Theorem \ref{para} and
Theorem \ref{bd-para}, we have $\sigma(U) = \sigma(V)$ and
$\sigma_{bd}(U) = \sigma_{bd}(V).$ Thus $\sigma(S) = \sigma(T)$
and  if $ \lambda \in \sigma_{bd}(T), $ then $ \frac{1}{\lambda}
\in \sigma_{bd}(U).$ Thus  $ \frac{1}{\lambda} \in \sigma_{bd}(V)$
and so $\lambda \in \sigma_{bd}(S).$ This implies that
$\sigma_{bd}(T) \subset \sigma_{bd}(S). $ Similarly we have
$\sigma_{bd}(S) \subset \sigma_{bd}(T) $ and then $\sigma_{bd}(T)=
\sigma_{bd}(S).$

\ep

 \baselineskip=12pt
\bigskip
\vspace{-5 mm }
 \baselineskip=12pt
\bigskip

{\tiny
\noindent Mohammed Berkani,\\
 \noindent Department of Mathematics,\\
 \noindent Science faculty of Oujda,\\
\noindent University Mohammed I,\\
\noindent Laboratory LAGA, \\
\noindent Morocco\\
\noindent berkanimo@aim.com,\\

\end{document}